\documentclass[10pt,twoside]{article}
\usepackage{graphicx}
\usepackage{amsmath}
\usepackage{Latex-document}
\usepackage{latexsym,psfig}
\usepackage{amsmath,amsfonts,amsthm,amssymb}
\newcommand{\R}{{\mathbb R}}
\renewcommand{\d}{\partial}

\newcommand{\e}{\epsilon}

\newcommand{\ty}{\tilde{y}}

\newcommand{\3}{\|\hspace{-.3ex}|}

\newtheorem{theorem}{Theorem}

 \newtheorem{definition}[theorem]{Definition}
  \newtheorem{conjecture}[theorem]{Conjecture}
 \newtheorem{open}[theorem]{Open Problem}

\title{\bf  Nonlinear Wave Equations\vskip 6mm}

\author{Daniel Tataru\vspace*{-0.5cm}\thanks{Department of Mathematics,
University of California at Berkeley, Berkeley, CA 94720, USA.
E-mail: tataru@math.berkeley.edu}}

\date{\vspace{-8mm}}

\begin{document}
\maketitle

\thispagestyle{first} \setcounter{page}{209}

\begin{abstract}\vskip 3mm
  The analysis of nonlinear wave equations has experienced a dramatic
  growth in the last ten years or so. The key factor in this has been
  the transition from linear analysis, first to the study of bilinear and
  multilinear wave interactions, useful in the analysis of semilinear
  equations, and next to the study of nonlinear wave interactions,
  arising in fully nonlinear equations.  The dispersion phenomena
  plays a crucial role in these problems. The purpose of this article
  is to highlight a few recent ideas and results, as well as to
  present some open problems and possible future directions in
  this field.

\vskip 4.5mm


\end{abstract}

\vskip 12mm

\section{Introduction} \label{section 1}\setzero
\vskip-5mm \hspace{5mm}

Consider the constant and variable coefficient wave
operators in $\R \times \R^n$,
\[
\Box = \d_t^2 - \Delta_x, \qquad \Box_g = g^{ij}(t,x) \d_i \d_j .
\]
In the variable coefficient case the summation occurs from $0$ to $n$
where the index $0$ stands for the time variable. To insure that the
equation is hyperbolic in time we assume that the matrix $g^{ij}$ has
signature $(1,n)$ and that the time level sets $t=const$ are
space-like, i.e.  $g^{00} > 0$. We consider semilinear wave equations,
\[
 \Box u = N(u) \qquad  (SLW), \qquad \Box u = N(u,\nabla u) \qquad (GSLW)
\]
and quasilinear wave equations,
\[
 \Box_{g(u)}  u = N(u)(\nabla u)^2 \qquad  (NLW), \qquad
 \Box_{g(u,\nabla u)}  u = N(u,\nabla u) \qquad (GNLW).
\]
To each of these equations we associate initial data in Sobolev spaces
\[
u(0) = u_0 \in H^s(\R^n), \qquad \d_t u(0) = u_1 \in
H^{s-1}(\R^n).
\]
There are two natural questions to ask: (i) Are the equations locally well-posed in $H^s \times H^{s-1}$?  (ii)
Are the solutions global, or is there blow-up in finite time?

\paragraph{Local well-posedness.}In a first approximation we
define it as follows:
\begin{definition}
  A nonlinear wave equation is well-posed in $H^s \times
  H^{s-1}$ if for each $(v_0,v_1) \in H^s \times H^{s-1}$ there is
  $T >0$ and a neighborhood $V$ of $(v_0,v_1)$ in $H^s \times
  H^{s-1}$ so that for each initial data $(u_0,u_1) \in V$ there is an
  unique solution $u \in C(-T,T;H^s)$, $\d_t u \in C(-T,T;H^{s-1})$
  which depends continuously on the initial data.
\end{definition}
In practice in order to prove uniqueness one often has to further
restrict the class of admissible solutions.  In most problems, the
bound $T$ from below for the life-span of the solutions can be chosen
to depend only on the size of the data.

It is not very difficult to prove that all of the above problems are
locally well-posed in $H^s \times H^{s-1}$ for large $s$. The
interesting question is what happens when $s$ is small. One indication
in this regard is given by scaling. At least in the case when the
nonlinear term has some homogeneity, for instance $N(u) = u^p$ or
$N(u) = u^p (\nabla u)^q$, one looks for an index $\alpha$ so that all
transformations of the form $u(x,t) \to \lambda^{\alpha} u(\lambda
x,\lambda t)$, $\lambda > 0$ leave the equation unchanged.
Correspondingly one finds an index $s_0 = \frac{n}{2} - \alpha$ so
that the norm of the initial data $(u_0,u_1)$ in the homogeneous
Sobolev spaces $\dot{H}^s \times \dot{H^{s-1}}$ is preserved by the
above transformations.

Below scaling ($s < s_0$) a small data small time result rescales into
a large data large time result. Heuristically one concludes that local
well-posedness should not hold. Still, to the author's knowledge there
is no proof of this yet.
\begin{conjecture}
Semilinear wave equations are ill-posed below scaling.
\end{conjecture}
This becomes much easier to prove if one strengthens the definition of
well-posedness, e.g.  by asking for uniformly continuous or $C^1$
dependence of the solution on the initial data.

If $s=s_0$ then for small initial data local well-posedness is
equivalent to global well-posedness. The same would happen for large
data if we were to strengthen the definition of well-posedness and ask
for a lifespan bound which depends only on the size of the data.
This is the only case where this distinction makes a difference.

If $s > s_0$ then a local well-posedness result gives bounds for
life-span $T_{max}$ of the solutions in terms of the size of the data,
\[
\| (u_0,u_1)\|_{H^s \times H^{s-1}} \leq M \implies T_{max}
\gtrsim M^{s_0-s}.
\]
The better localization in time makes the problems somewhat easier to
study.  However, besides scaling there are also other obstructions to
well-posedness.  These are related to various concentration phenomena
which can occur depending on the precise structure of the equation.

\paragraph{Global well-posedness.}

We briefly mention that there is a special case in which the global
well-posedness is well understood, namely when the initial data is
small, smooth, and decays at infinity. This is not discussed at all in
what follows.

Consider first the case when $s$ is above scaling, $s > s_0$, and
local well-posedness holds in $H^s \times H^{s-1}$.  Then any solution
can be continued as long as its size does not blow-up. Hence the goal
of any global argument should be to establish a-priori bounds on the
$H^s \times H^{s-1}$ norm of the solution. All known results of this
type are for problems for which there are either conserved or
quasi-conserved positive definite quantities. Such conserved
quantities can often be found for equations which are physically
motivated or which have some variational structure. For simplicity
suppose that there is some index $s_c$ and an energy functional $E$ in
$H^{s_c} \times H^{s_c -1}$ which is preserved along the flow. The
index $s_c$ needs not be equal to the scaling index $s_0$.
There are three cases to consider:

(i) The subcritical case $s_c > s_0$. Then a local well-posedness result
at $s = s_c$ implies the global result for $s \geq s_c$. Furthermore,
in recent years there has been considerable interest
in establishing global well-posedness also for $s_0 < s < s_c$.
This is based on an idea first introduced by Bourgain~\cite{Bo1}
in a related problem for the Schr\"oedinger equation, and followed up
by a number of authors.

(ii) The critical case $s_c=s_0$. Here the energy is not needed
for small data, when local and global well-posedness are equivalent.
For large data, however, the energy conservation is not sufficient
in order to establish the existence of global solutions. In addition,
one needs a non-concentration argument, which should say that
the energy cannot concentrate inside a characteristic cone.

(iii) The supercritical case, $s_c  < s_0$. No  global results are
known:
\begin{open}
{\rm Are supercritical problems globally well-posed for $s \geq s_0$?}
\end{open}
A simple example is the equation (NLW) with $N(u)=|u|^{p-1} u$.  The energy is
\[
E(u) = \int |u_t|^2+|\nabla_x u|^2 +\frac{1}{p+1} |u|^{p+1}.
\]
Then $s_c = 1$, while $s_0 = \frac{n}{2} - \frac{2}{p-1}$. In $3+1$
dimensions, for instance, $p=3$ is subcritical, therefore one has
global well-posedness in $H^1\times L^2$. The exponent $p=5$ is
critical and  in this case the problem is known to be globally well-posed
in $ H^1\times L^2$; the non-concentration argument is due to
Grillakis~\cite{G}. The exponent $p=7$ is supercritical.

\paragraph{Blow-up.}
Not all nonlinear wave equations are expected to have
global solutions. Quite the contrary, generic equations are expected
to blow up in finite time; only for problems with some special
structure it seems plausible that global well-posedness may hold.  A
simple way to produce blow-up is to look for self-similar solutions,
$u(x,t) = t^\gamma u(\frac{x}{t})$. If they exist, self-similar
solutions disprove global well-posedness. Because they must respect
the scaling of the problem, they are not so useful when trying to
disprove local well-posedness.

Another way to produce blow-up solutions is the so-called ode blow-up.
In the simplest setting this means looking at one dimensional
solutions (say $u(x,t)=u(t)$) which solve an ode and blow up in finite
time.  Then one can truncate the initial data spatially and still
retain the blow-up because of the finite speed of propagation. This is
still not very useful for the local problem.

A better idea is to constructs blow-up solutions which are
concentrated essentially along a light ray, see
Lindblad~\cite{L0},\cite{L1} and Alihnac~\cite{A}.  In this setup the
actual blow-up occurs either because of the increase in the amplitude,
in the semilinear case, or because of the focusing of the light rays,
in the quasilinear case. As it turns out, the counterexamples of this
type are often sharp for the local well-posedness problem.

\section{Semilinear wave equations} \label{section 2}
\setzero\vskip-5mm \hspace{5mm}

Usually, a fixed point argument is used to obtain local results for
semilinear equations. We first explain this for  the case when
$s=s_0$. We define the homogeneous and inhomogeneous solution
operators, $S$ and $\Box^{-1}$ by
\[
S(u_0,u_1) = u \Longleftrightarrow \{ \Box u = 0, \quad u(0) =
u_0, \quad \d_t u(0) = u_1  \},
\]
\[
\Box^{-1}f = u \Longleftrightarrow \{ \Box u = f, \quad u(0) = 0,
\quad \d_t u(0) = 0  \}.
\]
Then the equation (NLW) for instance can be recast as
\[
u = S(u_0,u_1) + \Box^{-1} N(u).
\]
To solve this using a fixed point argument one needs two Banach spaces
$X$ and $Y$ with the correct scaling and the following mapping
properties:
\[
S:H^s \times H^{s-1} \to X, \quad \Box^{-1}: Y \to X, \quad
N:X\to Y.
\]
The first two are linear, but the last one is nonlinear.  The small
Lipschitz constant is always easy to obtain provided the initial data
is small and that $N$ decays faster than linear at $0$. The solutions
given by the fixed point argument are global.

In the case $s > s_0$ the scaling is lost, and with this method one can only hope to
get  results which are local in time. To localize in time one
chooses a smooth compactly supported cutoff function $\chi$ which
equals $1$ near the origin. The fixed point argument is now used for
the equation
\[
u = \chi S(u_0,u_1) + \chi \Box^{-1} N(u).
\]
A solution to this solves the original equation only in an interval
near the origin where $\chi = 1$. The modified mapping properties
are
\[
\chi S:H^s \times H^{s-1} \to X, \quad \chi \Box^{-1}: Y \to X,
\quad N:X\to Y.
\]
How does one choose the spaces $X$, $Y$? One approach is to use the energy estimates for the wave equation and set
\[
X = \{ u \in L^\infty(H^s), \nabla u \in L^\infty(H^{s-1})\},
\quad Y = L^1(H^{s-1}).
\]
The first two mapping properties are trivial. However, if
the third holds then we must also have $N: X \to  L^\infty(H^{s-1})$.
The one unit difference in scaling between $L^1$ and $L^\infty$
implies that this can only work for $s \geq s_0+1$.

What is neglected in the above setup is the dispersive properties of
the wave equation. Solutions to the linear wave equation cannot stay
concentrated for long time intervals. Instead, they will disperse and
decay in time (even though the energy is preserved). In harmonic
analysis terms, this is related to the restriction theorem (see
\cite{S2}) and is a consequence of the nonvanishing curvature of the
characteristic set for the wave operator, namely the cone $\xi_0^2 =
\xi_1^2 + \cdots +\xi_n^2$. Here $\xi$ stands for the Fourier
variable. One way of quantifying the dispersive effects is through the
Strichartz estimates. They apply both to the homogeneous and the
inhomogeneous equation (see \cite{KT} and references therein):
\[
S: H^{\rho} \times H^{\rho -1} \to L^pL^q, \qquad
|D|^{1-\rho_1-\rho} \Box^{-1}: L^{p'_1} L^{q'_1} \to L^{p}L^{q}
\]
where $(\rho,p,q)$ and $(\rho_1,p_1,q_1)$ are subject to
\[
\frac{1}{p}+\frac{n}{q} = \frac{n}{2}-\rho, \quad
\frac{2}{p}+\frac{n-1}{q} \leq \frac{n-1}{2}, \quad 2 \leq p,q
\leq \infty, \quad (\rho,p,q) \neq (1,2,\infty).
\]
The worst case in these estimates occurs for certain highly localized
approximate solutions to the wave equation, which are called wave
packets.  A frequency $\lambda$ wave packet on the unit time scale is
essentially a bump function in a parallelepiped of size $1 \times
\lambda^{-1} \times (\lambda^{-\frac12})^{n-1}$ which is obtained from
a $\lambda^{-1} \times (\lambda^{-\frac12})^{n-1}$ parallelepiped at
time zero which travels with speed $1$ in the normal direction.
Because of the uncertainty principle, this is the best possible
spatial localization which remains coherent up to time $1$. Of course
one can rescale and produce wave packets on all time scales.

In low dimension $n=2,3$ the Strichartz estimates provide a complete
set of results for generic equations of both (NLW) and (GNLW) type.
Consider the following two examples, of which the second is wrong but
almost right:
\[
\Box u = u^3, \quad n=3, \quad s=s_0=\frac12, \quad X=L^4, \quad
Y =L^{\frac43},
\]
\[
\Box u = u \nabla u, \quad n=3, \quad s_0=\frac12, \quad s = 1
\quad X=|D|^{-1}L^\infty L^2 \cap  L^2 L^\infty \quad Y =L^2 .
\]
For $n \geq 4$, however, the Strichartz estimates no longer provide
all the results. The reason is as follows. The worst nonlinear
interaction in both (NLW) and (GNLW) occurs for wave packets which
travel in the same direction.  One can use the Strichartz estimates to
accurately describe the interaction of same frequency wave packets.
But in the interaction of two wave packets at different frequencies,
the low frequency packet is more spread, and only a small portion of
it will interact with the high frequency packet.  However, unlike in
low dimension, the Strichartz estimates do not provide sharp bounds
for this smaller part of a wave packet.

A more robust idea due to Bourgain~\cite{Bo0} and
Klainerman-Machedon~\cite{KM} is to use the $X^{s,b}$ spaces
associated to the wave equation very much in the same way the Sobolev
spaces are associated to the Laplacian:
\[
\| u\|_{X^{s,b}} = \| (1+|\xi|)^s (1+||\xi_0| - |\xi'||)^b
\hat{u}\|_{L^2} .
\]
Then one chooses $X=X^{s,\frac12}$ and $Y=X^{s-1,-\frac12}$.  The
Strichartz information is not lost since for $\rho,p,q$ as above we
have the dual embeddings
\[
X^{\rho,\frac12+} \subset L^pL^q, \qquad L^{p'}L^{q'} \subset
X^{-\rho,-\frac12-}.
\]
Within the framework of the $X^{s,b}$ spaces one can prove bilinear
estimates which provide a better description of the interaction of
high and low frequencies, see \cite{MR2001g:35145} and references
therein.  The bilinear estimates are obtained as weighted convolution
estimates in the Fourier space, by using the above embeddings, or by
combining the two methods.  Sometimes even this setup does not suffice
and has to be modified further, see \cite{MR1739207}.
\begin{conjecture}
The equation $\Box u = u^p$ is locally well-posed in $H^s \times
H^{s-1}$ for $n \geq 4$, $0 \leq  s \leq \frac12$, $p
(\frac{n+1}{4}-s) \leq  (\frac{n+5}{4}-s)$. {\rm (see
\cite{MR2000a:35175}
 for more details)}
\end{conjecture}
\paragraph{The null condition.}
A natural question to ask is whether there are equations which behave
better than generic ones. This may happen if the worst interaction
(between parallel wave packets) does not occur in the nonlinearity.
A good example is (GNLW) with a quadratic nonlinearity $Q(\nabla u,
\nabla u)= q^{ij} \d_i u\, \d_j u$. The cancellation condition, called
null condition, asserts that
\[
q^{ij} \xi_i \xi_j = 0 \qquad \text{in the characteristic set} \
g^{ij} \xi_i \xi_j = 0.
\]
All such null forms are linear combinations of
\[
Q_{ij}(\nabla u,\nabla v) = \d_i u \d_j v- \d_i v \d_j u, \qquad
Q_0(u,v)= g^{ij} \d_i u\d_j v.
\]
\begin{open}
{\rm Study semilinear wave equations corresponding to variable
coefficient wave operators for $n \geq 4$ (generic case) or $n
\geq 2$ (with null condition).}
\end{open}

In the constant coefficient case one can easily use the null condition
in the context of the $X^{s,b}$ spaces.  This is done using
inequalities of the following form:
\[
|q_0(\xi,\eta)| \leq c(|p(\xi)|+|p(\eta)|+|p(\xi+\eta)|)
\]
respectively
\[
|q_{ij}(\xi,\eta)| \leq c|\xi|^\frac12|\eta|^{\frac12}|\xi+\eta|^\frac12
(|p(\xi)|^{\frac12} + |p(\eta)|^\frac12 + |p(\xi+\eta)|^\frac12)
\]
where by $p(\xi)$ we denote the symbol of the constant coefficient wave
operator, given by $p(\xi) = \xi_0^2-\xi_1^2-\cdots -\xi_n^2$.
Combining this with the embeddings above one can lower the $s$ in the
local theory whenever the null condition is satisfied.  Unfortunately,
this does not always give optimal results.  The problem of obtaining
improved $L^pL^q$ estimates for null forms has also been explored, see
\cite{W}\cite{MR1865417} \cite{bilinear}, but without immediate
applications to semilinear wave equations.  We limit the following
discussion to two of the more interesting models.

\paragraph{Wave maps.} These are functions from $\R^n\times \R$
into a complete Riemannian manifold $(M,g)$ which are critical points for
\[
I (\phi) = \int_{\R^n\times \R} |\d_t \phi|_g^2 - |\nabla_x
\phi|_g^2 dx \ dt .
\]
In local coordinates the equation for wave maps has the form
\[
\Box \phi^k = \Gamma^{k}_{ij}(\phi) Q_0(\phi^i \phi^j)
\]
where $\Gamma^{k}_{ij}$ are the Riemann-Christoffel symbols.
The energy functional is
\[
E(u) =  \int_{\R^n} |\d_t \phi|_g^2 + |\nabla_x \phi|_g^2 dx.
\]
The scaling index is $s_0 = \frac{n}{2}$ and $s_c = 1$.  Local
well-posedness for $s > s_c$ can be obtained using the $X^{s,b}$
spaces. For $s = s_c$, using some modified $X^{s,b}$ spaces, local
(and therefore small data global) well-posedness was established first
in homogeneous Besov spaces $B^{\frac{n}2}_{2,1} \times
B^{\frac{n}2-1}_{2,1}$ in Tataru~\cite{2002c:58045} and then in
Sobolev spaces by Tao~\cite{MR1869874} (for the sphere, $n \geq 2$)
and other authors (general target manifold, $n \geq 3$).  Large data
global well-posedness is false in the supercritical case $n \geq 3$,
where self-similar blowup can occur. This leaves open problems in the
critical case $n=2$:

\begin{conjecture}
(i) The two dimensional wave maps equation is globally  well-posed for
small data in $H^1 \times L^2$ for any complete target manifold.

(ii) The two dimensional wave maps equation is globally well-posed for
large data in $H^1 \times L^2$ for ``good'' target manifolds.
\end{conjecture}

\paragraph{ The Yang Mills equations.}
Given a compact Lie group $\cal G$ whose Lie algebra $\mathbf g$
admits an invariant inner product $\langle \cdot,\cdot \rangle$ one
considers $\mathbf g$ valued connection $1$-forms $A_j d x^j$ in
$\R^n\times \R$.  The covariant derivatives of $\mathbf g$ valued
functions are defined by
\[
D_j B = \d_j B + [A_j,B].
\]
The ($\mathbf g$ valued) curvature of the connection $A$ is
\[
F_{ij} = \d_i A_j-\d_j A_i + [A_i, A_j].
\]
This is invariant with respect to gauge transformations
\[
A_j \to O A_j O^{-1} - \d_j OO^{-1}, \qquad O \in {\cal G}.
\]
A Yang-Mills connection is a critical point for the Yang-Mills
functional
\[
I (A) = \int_{\R^n\times \R} \langle F_{ij}, F^{ij} \rangle
dx \ dt
\]
where indices are lifted with respect to the Minkovski metric.
Then the Yang-Mills equations have the form
\[
D^j F_{ij} = 0
\]
and the energy functional is
\[
E(A) = \int_{\R^n\times \R} \langle F_{ij}, F_{ij} \rangle dx \
dt.
\]
A Yang-Mills connection is not a single connection, but instead it is
a class of equivalence with respect to the above gauge transformation.
In order to view the Yang-Mills equations as semilinear wave equations
and solve them one has to fix the gauge, i.e. select a single
representative out of each equivalence class. Common gauge choices
include: (i) the temporal gauge $A_0=0$, (ii) the wave gauge
$\d_j A^j = 0$ and (iii) the Coulomb gauge $\sum_{j=1}^n \d_j A_j = 0$.
To understand the equation better it may help
to look first at an  oversimplified version, namely
\[
\Box u = (u \cdot \nabla_x)u + \nabla_x p, \qquad \nabla_x \cdot
u = 0.
\]
This exhibits a $Q_{ij}$ type null condition. The scaling index is
$s_0 =\frac{n-2}{2}$ and $s_c=1$.  Using the $X^{s,b}$ spaces one can
improve the local theory somewhat, but certain more subtle
modifications of this are needed in order to handle high-low frequency
interactions. In \cite{MR1626261} such an approach is used
to prove that local well-posedness holds for $s > s_0$, $n \geq 4$.

\begin{open}
{\rm Is the Yang-Mills equation well-posed for $s > s_0$, $n=2,3$? (Likely not for $n=2$. For $n=3$ one can obtain
$s > \frac34$ using the $X^{s,b}$ spaces.)}
\end{open}

\begin{conjecture}
  (i) The Yang-Mills equation is globally well-posed for small data in
  $H^{s_0} \times H^{s_0-1}$ for $n \geq 4$.

  (ii) The Yang-Mills equation is globally well-posed for large data
  in $H^{s_0} \times H^{s_0-1}$ for $n = 4$.
\end{conjecture}

\section{Nonlinear wave equations} \label{section 3} \setzero\vskip-5mm \hspace{5mm}

We consider (NLW), since (GNLW) reduces to it by differentiation. The fixed
point argument in the semilinear case cannot be applied in the
nonlinear case, because the wave equation parametrix is not strongly
stable with respect to small changes in the coefficients. Instead, one
must adopt a different strategy: (i) show that local solutions exist
for smooth data, (ii) obtain a-priori bounds for smooth solutions
uniformly with respect to initial data in a bounded set in $H^s \times
H^{s-1}$ and (iii) prove continuous dependence on the data in a weaker
topology, and obtain solutions for $H^s \times H^{s-1}$ data as weak
limits of smooth solutions.  Steps (i) and (iii) are more or less
routine, it is (ii) which causes most difficulties. A good starting
point is Klainerman's energy estimate \vspace*{-2ex}
\[
\|\nabla u(t)\|_{H^{s-1}} \lesssim \|\nabla u(0)\|_{H^{s-1}} \exp{( \int_0^t \|\nabla
  u(s)\|_{L^\infty}}ds). \vspace*{-2ex}
\]
This shows that all Sobolev norms of a solution remain bounded for as
long as $\|\nabla u\|_{L^1L^\infty}$ stays bounded.  It remains to see
how to obtain bounds on $\|\nabla u\|_{L^1L^\infty}$. The
classical approach uses energy estimates and Sobolev embeddings,
but, as in the semilinear case, it only yields results one
unit above scaling, namely for $s > \frac{n}{2} +1$.

Better results could be obtained using the Strichartz estimates
instead. However, this is very nontrivial as one would have to
establish the Strichartz estimates for the operator $\Box_{g(u)}$, which
has very rough coefficients. Compounding the difficulty,
the argument is necessarily circular,
\smallskip

\parbox{.8in}{coefficients \\ regularity} $\implies$ \parbox{.8in}{Strichartz \\
  estimates} $\implies$
\parbox{.8in}{solution\\ regularity} $\implies$ \parbox{.8in}{coefficients\\  regularity}
\smallskip

One can get around this with a bootstrap argument of the form
\[
\left.\begin{array}{c} \|(u_0,u_1)\|_{H^s \times H^{s-1}} \leq \e \cr
\3 g(u) \3 \leq 2 \end{array} \right\} \implies \left\{\begin{array}{c}
\text{Strichartz estimates for $\Box_{g(u)}$ in $[-1,1]$}
\cr \3 g(u) \3 \leq 1
\end{array}\right.
\]
where the (possibly nonlinear) triple norm contains the needed
information about the metric. Still, apriori there is no clear way to
determine exactly how it should be defined. A starting point is to set
$\3 g(u) \3 = \|\nabla g\|_{L^1L^\infty}$, but this only leads to
partial results.  Following partial results independently obtained by
Bahouri-Chemin~\cite{BC},\cite{BC1} and Tataru~\cite{nlw},~\cite{lp}
and further work of Klainerman-Rodnianski~\cite{KR}, the next result
represents the current state of the problem:

\begin{theorem} {\rm (Smith-Tataru~\cite{ST})} The equation (NLW) is locally
  well-posed in $H^{s} \times H^{s-1}$ for $s > \frac{n}{2} +\frac34$
  ($n=2$) and $s > \frac{n}{2} +\frac12$ ($n=3,4,5$). In addition,
the Strichartz estimates with $q=\infty$ hold for the corresponding
wave operator $\Box_{g(u)}$.
\end{theorem}
Lindblad's counterexamples correspond to $s = \frac{n+3}{4}$ and show
that this result is sharp for $n=2,3$.  The restriction to $n \leq 5$
is not central to the problem, it can likely be removed with some
extra work.
\begin{open}
{\rm Improve the above result in dimension $n \geq 4$.}
\end{open}
\paragraph{ Wave equation parametrices.}
In most approaches, the key element in the proof of the Strichartz
estimates is the construction of a parametrix for the wave equation.
There are many ways to do this for smooth coefficients, however, as
the regularity of the coefficients decreases, they start to break
down. Let us begin with the classical Fourier integral operator
parametrix, used in the work of  Bahouri-Chemin:
\[
K(x,y) = \int a(x,y,\xi) e^{i \phi(x,y,\xi)} d \xi .
\]
The phase $\phi$ is initialized by $\phi(x,y,\xi) = \xi (x-y)$ when
${x_0=y_0}$ and must solve an eikonal equation, while for the
amplitude $a$ one obtains a transport equation along the Hamilton
flow. The disadvantage is that all spatial localization comes from
stationary phase, which seems to require too much regularity for the
coefficients.

One way to address the issue of spatial localization is to begin with
wave packets, which have the best possible spatial localization on the
unit time scale. In the variable coefficient case the frequency
$\lambda$ wave packets are bump functions on curved parallelepipeds of
size $1 \times \lambda^{-1} \times (\lambda^{-\frac12})^{n-1}$. These
parallelepipeds are images of $ \lambda^{-1} \times
(\lambda^{-\frac12})^{n-1}$ parallelepipeds at the initial time,
transported along the Hamilton flow for $\Box_g$ corresponding to
their conormal direction. Then one can seek approximate solutions for
$\Box_g$ as discrete superpositions of wave packets, $u = \sum_T u_T$.
It is not too difficult to construct individual wave packets, the more
delicate point is to show that the wave packets are almost orthogonal.
This approach, which is used in \cite{ST}, was originally introduced
by Smith~\cite{Sm} and used to prove the Strichartz estimates in $2$
and $3$ dimensions for operators with $C^2$ coefficients.

Another parametrix with a better built in spatial localization can be
obtained by doing a smooth phase space analysis:
\[
K(\ty,y) = \int_C a(x,\xi) e^{i
(\phi(y,x,\xi)-\overline{\phi(\ty,x_t,\xi_t)}} dx\ d \xi\ dt
\quad \phi(y,x,\xi) = \xi(x-y) + i|\xi| (x-y)^2.
\]
Here $(x,\xi) \to (x_t,\xi_t)$ is the Hamilton flow for $\Box_g$ on
the characteristic cone $C= \{g^{ij}(x) \xi_i \xi_j = 0\}$.  One can
factor this into a product of three operators, namely an FBI
transform, a phase space transport along the Hamilton flow and then an
inverse FBI transform. Neglecting the first one, i.e. setting $x=y$
above, produces an operator which is similar to the Fourier integral
operators with complex phase. However, it seems to be more useful to
keep the Gaussian localizations at both ends. Parametrices of this
type were introduced in Tataru~\cite{cs} and used to prove Strichartz
estimates for operators with $C^2$ coefficients in all dimensions. The
$C^2$ condition was later relaxed in \cite{lp} to $\nabla^2 g \in
L^1L^\infty$. Localization and scaling arguments lead also to weaker
estimates for operators whose coefficients have less regularity.  Such
estimates are known to be sharp, see the counterexamples in
Smith-Tataru~\cite{STcounter}.

\paragraph{The null condition.}

As in the semilinear case, one may ask whether better results can be
obtained for equations with special structure. However, unlike the
semilinear case, little is known so far. We propose the following

\begin{definition}
We say that the equation (GNLW) satisfies the null condition if
\[
\frac{\d g^{ij}(u,p)}{\d p_k} \xi_i \xi_j \xi_k = 0 \qquad
\text{in} \ g^{ij}(u,p) \xi_i \xi_j = 0.
\]
\end{definition}

\begin{conjecture}
If the null condition holds then the equation (GNLW) is well-posed
in $H^s \times H^{s-1}$ for some $s < \frac{n}2 + \frac34$ ($n=2$)
respectively  for some $s <  \frac{n}2 + \frac12$ ($n = 3$).
\end{conjecture}

In $3+1$ dimensions a problem which does not quite fit into the above
setup but still satisfies some sort of null condition is the
Einstein's equations in general relativity. It is similar to the Yang
Mills equations in that it has a gauge invariance, and the null
condition is only apparent after fixing the gauge.
Klainerman-Rodnianski have obtained a different proof of
Theorem 9 for this special case of (NLW).

\label{lastpage}

\end{document}